\documentclass[twoside,11pt]{article}
\usepackage{latexsym,amssymb,times} \input amssym.def \input amssym
\newtheorem{te}{Theorem}

\newtheorem{pro}{Proposition}
\newtheorem{lem}{Lemma}
\newtheorem{df}{Definition}

\newtheorem{pb}{Problem}
\def\bt{\hspace{-2mm}{\bf .}\hspace{2mm}}
\def\dok{\noindent{\bf Proof. }}
\def\kdok{\hfill $\Box$ \par \vspace*{2mm} }
\def\dom{\mathop{\rm dom}\nolimits}

\begin{document}
\thispagestyle{plain}
\begin{center}
           {\large \bf FOUR GAMES ON BOOLEAN ALGEBRAS }
\end{center}
\begin{center}
{\small \bf Milo\v s S.\ Kurili\'c and Boris  \v Sobot}\\[2mm]
{\small  Department of Mathematics and Informatics, University of Novi Sad,\\
Trg Dositeja Obradovi\'ca 4, 21000 Novi Sad, Serbia\\
e-mail: milos@dmi.uns.ac.rs, sobot@dmi.uns.ac.rs}
\end{center}
\begin{abstract}
\noindent The games ${\cal G}_2$ and ${\cal G}_3$ are played on a complete Boolean algebra $\Bbb B$ in
$\omega$-many moves. At the beginning White picks a non-zero element $p$ of $\Bbb B$ and, in the $n$-th move,
White picks a positive $p_n<p$ and Black chooses an $i_n\in\{0,1\}$. White wins ${\cal G}_2$ iff $\liminf
p_n^{i_n}=0$ and wins ${\cal G}_3$ iff $\bigvee_{A\in [\omega ]^{\omega }}\bigwedge_{n\in A}p_n^{i_n}=0$. It is
shown that White has a winning strategy in the game ${\cal G}_2$ iff White has a winning strategy in the
cut-and-choose game ${\cal G}_{\rm{c\&c}}$ introduced by Jech. Also, White has a winning strategy in the game
${\cal G}_3$ iff forcing by $\Bbb B$ produces a subset $R$ of the tree ${}^{< \omega }2$ containing either
$\varphi ^{\smallfrown } 0 $ or $\varphi ^{\smallfrown }1$, for each $\varphi\in {}^{<\omega }2$, and having
unsupported intersection with each branch of the tree ${}^{< \omega }2$ belonging to $V$. On the other
hand, if forcing by $\Bbb B$ produces independent (splitting) reals then White has a winning strategy in the game ${\cal
G}_3$ played on $\Bbb B$. It is shown that $\diamondsuit$ implies the existence of an algebra on which these
games are undetermined.
\vspace{1mm}
\end{abstract}

\section*{{\normalsize \bf 1. Introduction}}

In \cite{Jech2} Jech introduced the cut-and-choose game ${\cal G}_{\rm{c\&c}}$, played by two players, White
and Black, in $\omega$-many moves on a complete Boolean algebra ${\Bbb{B}}$ in the following way. At the
beginning, White picks a non-zero element $p\in{\Bbb{B}}$ and, in the $n$-th move, White picks a non-zero
element $p_n<p$ and Black chooses an $i_n\in\{0,1\}$. In this way two players build a sequence $\langle
p,p_0,i_0,p_1,i_1,\dots\rangle$ and White wins iff $\bigwedge_{n\in\omega }p_n^{i_n}=0$ (see Definition
\ref{D900}).

A winning strategy for a player, for example White,
is a function which, on the basis of the previous moves of both
players, provides ``good" moves for White such that White always wins.
So, for a complete Boolean algebra ${\Bbb{B}}$ there are three possibilities:
1) White has a winning strategy; 2) Black has a winning strategy or
3) none of the players has a winning strategy.
In the third case the game is said to be undetermined on ${\Bbb{B}}$.

The game-theoretic properties of Boolean algebras have interesting algebraic and forcing translations. For
example, according to \cite{Jech2} and well-known facts concerning infinite distributive laws we have the
following results.

\begin{te}\bt\rm\label{T914}
(Jech)
For a complete Boolean algebra ${\Bbb{B}}$ the following conditions are equivalent:

(a) White has a winning strategy in the game ${\cal G}_{\rm{c\&c}}$;

(b) The algebra ${\Bbb{B}}$ does not satisfy the $(\omega ,2)$-distributive law;

(c) Forcing by ${\Bbb{B}}$ produces new reals in some generic extension;

(d) There is a countable family of 2-partitions of the unity having no common refinement.
\end{te}

Also, Jech investigated the existence of a winning strategy for Black and using $\diamondsuit $
constructed a Suslin algebra in which the game ${\cal G}_{\rm{c\&c}}$ is undetermined.
Moreover in \cite{Zapl1} Zapletal  gave a
ZFC example of a complete Boolean algebra in which the game ${\cal G}_{\rm{c\&c}}$  is undetermined.

Several generalizations of the game ${\cal G}_{\rm{c\&c}}$ were considered. Firstly, instead of cutting of $p$
into two pieces, White can cut into $\lambda $ pieces and Black can choose more than one piece (see \cite{Jech2}).
Secondly, the game can be of uncountable length so Dobrinen in \cite{Dobr1} and \cite{Dobr2} investigated the
game ${\cal G} _{<\mu}^{\kappa }(\lambda )$ played in $\kappa $-many steps in which White cuts into $\lambda $ pieces
and Black chooses less then $\mu$ of them.

In this paper we consider three games ${\cal G}_2 , {\cal G}_3$ and ${\cal G}_4$  obtained from the game
${\cal G}_{\rm{c\&c}}$ (here denoted by ${\cal G}_1$) by changing the winning criterion in the following way.

\begin{df} \bt \rm   \label{D900}
The games ${\cal G}_k $, $k\in \{ 1,2,3,4\} $,
are played by two players, White and Black,
on a complete Boolean algebra ${\Bbb{B}}$ in $\omega $-many moves. At the beginning White chooses a
non-zero element $p\in{\Bbb{B}}$. In the $n$-th move White chooses a $p_n\in(0,p)_{{\Bbb{B}}}$
and Black responds choosing $p_n$ or $p\setminus p_n$ or, equivalently,
picking an $i_n\in\{0,1\}$ chooses
$p_n^{i_n}$, where, by definition, $p_n^0=p_n$ and $p_n^1=p\setminus p_n$.
White wins the play  $\langle p,p_0,i_0,p_1,i_1,\dots\rangle$ in the game

\vspace{2mm}

${\cal G}_1 $  if and only if $\bigwedge_{n\in\omega } p_n^{i_n}=0$;

${\cal G}_2 $  if and only if $\bigvee_{k\in\omega }\bigwedge_{n\geq k}p_n^{i_n}=0$,
                              that is $\liminf p_n^{i_n}=0$;

${\cal G}_3 $  if and only if $\bigvee_{A\in [\omega ]^{\omega }}\bigwedge_{n\in A}p_n^{i_n}=0$;

${\cal G}_4 $  if and only if $\bigwedge_{k\in\omega }\bigvee_{n\geq k}p_n^{i_n}=0$,
                              that is $\limsup p_n^{i_n}=0$.

\end{df}
In the following theorem we list some results concerning the game ${\cal  G}_4$
which are contained in \cite{KuSo}.

\begin{te}\bt\rm\label{T910}
(a) White has a winning strategy in the game ${\cal  G}_4    $ played on a complete Boolean algebra ${\Bbb{B}}$
iff forcing
by ${\Bbb{B}}$ collapses $\goth c$ to $\omega$ in some generic extension.

(b) If ${\Bbb{B}}$ is the Cohen algebra $\rm{r.o.}({}^{<\omega}2,\supseteq)$ or a Maharam algebra (i.e. carries
a positive Maharam submeasure) then Black has a winning strategy in the game ${\cal  G}_4    $ played on ${\Bbb{B}}$.

(c) $\diamondsuit$ implies the existence of a Suslin algebra on which the game ${\cal  G}_4    $ is undetermined.
\end{te}

The aim of the paper is to investigate the game-theoretic properties of complete Boolean algebras related to the games
${\cal G}_2$ and ${\cal G}_3$. So, Section 2 contains some technical results,
in Section 3 we consider the game ${\cal G}_2$,
Section 4 is devoted to the game ${\cal G}_3$ and Section 5 to the algebras
on which these games are undetermined.

Our notation is standard and follows \cite{Jech3}. A subset of $\omega$ belonging to a generic extension will be called
supported iff it contains an infinite subset of $\omega$ belonging to the ground model. In particular, finite subsets of $\omega$ are unsupported.

\section*{{\normalsize \bf 2. Winning a play, winning all plays}}

Using the elementary properties of Boolean values and forcing it is easy to prove the following two statements.

\begin{lem}\bt\rm\label{T905}
Let ${\Bbb{B}}$ be a complete Boolean algebra, $\langle b_n:n\in\omega\rangle$ a sequence in ${\Bbb{B}}$ and
$\sigma=\{\langle\check n,b_n\rangle:n\in\omega\}$ the corresponding name for a subset of $\omega$. Then

(a) $\bigwedge_{n\in\omega}b_n=\|\sigma=\check\omega\|$;

(b) $\liminf b_n=\|\sigma\mbox{ is cofinite}\|$;

(c) $\bigvee_{A\in[\omega]^{\omega}}\bigwedge_{n\in A}b_n=\|\sigma\mbox{ is supported}\|$;

(d) $\limsup b_n=\|\sigma\mbox{ is infinite}\|$.
\end{lem}

\begin{lem}\bt\rm\label{T906}
Let ${\Bbb{B}}$ be a complete Boolean algebra, $p\in{\Bbb{B}}^+$, $\langle p_n:n\in\omega\rangle$ a sequence in
$(0,p)_{{\Bbb{B}}}$ and $\langle i_n:n\in\omega\rangle\in{}^{\omega}2$. For $k\in\{0,1\}$ let
$S_k=\{n\in\omega:i_n=k\}$ and let the names $\tau$ and $\sigma$ be defined by $\tau=\{\langle\check
n,p_n\rangle:n\in\omega\}$ and $\sigma=\{\langle\check n,p_n^{i_n}\rangle:n\in\omega\}$. Then

(a) $p'\Vdash\tau=\sigma= \check\emptyset$;

(b) $p\Vdash\tau=\sigma\triangle\check S_1$;

(c) $p\Vdash\sigma=\tau\triangle\check S_1$;

(d) $p\Vdash\sigma=\check\omega\Leftrightarrow\tau=\check S_0$;

(e) $p\Vdash\sigma=^*\check\omega\Leftrightarrow\tau=^*\check S_0$;

(f) $p\Vdash |\sigma|<\check\omega\Leftrightarrow\tau=^*\check S_1$.
\end{lem}

\begin{te}\bt\rm\label{T907}
Under the assumptions of Lemma \ref{T906}, White wins the play $\langle p,p_0,i_0,$ $p_1,i_1,\dots\rangle$ in
the game

${\cal G}_1$ iff $\| \sigma\mbox{ is not equal to }\check\omega \| =1$ iff $p\Vdash\tau\neq \check S_0$;

${\cal G}_2$ iff $\| \sigma \mbox{ is not cofinite}\| =1$ iff $p\Vdash\tau\neq^* \check S_0$;

${\cal G}_3$ iff $\| \sigma \mbox{ is not supported}\| =1$ iff $p\Vdash``\tau\cap \check S_0 \mbox{ and }\check
S_1\setminus\tau$ are unsupported";

${\cal G}_4$ iff $\| \sigma \mbox{ is not infinite}\| =1$ iff $p\Vdash\tau=^* \check S_1$.
\end{te}

\dok
We will prove the statement concerning the game ${\cal G}_3$ and leave the rest to the reader. So, White wins
${\cal G}_3$ iff $\bigvee_{A\in[\omega]^{\omega}}\bigwedge_{n\in A}p_n^{i_n}=0$, that is, by Lemma \ref{T905},
$\|\sigma\mbox{ is not supported}\|=1$ and the first equivalence is proved.

Let $1\Vdash``\sigma$ is not supported" and let $G$ be a ${\Bbb{B}}$-generic filter over $V$ containing
$p$. Suppose $\tau_G\cap S_0$ or $S_1\setminus\tau_G$ contains a subset
$A\in[\omega]^{\omega}\cap V$. Then $A\subseteq\sigma_G$, which is impossible.

On the other hand, let $p\Vdash``\tau\cap\check S_0\mbox{ and }\check S_1\setminus\tau$ are unsupported" and let
$G$ be a ${\Bbb{B}}$-generic filter over $V$. If $p'\in G$ then, by Lemma \ref{T906}(a), $\sigma_G=\emptyset$ so
$\sigma_G$ is unsupported. Otherwise $p\in G$ and by the assumption the sets $\tau_G\cap S_0$ and
$S_1\setminus\tau_G$ are unsupported. Suppose $A\subseteq\sigma_G$ for some $A\in[\omega]^{\omega}\cap V$. Then
$A=A_0\cup A_1$, where $A_0=A\cap S_0\cap\tau_G$ and $A_1=A\cap S_1\setminus\tau_G$, and at least one of these
sets is infinite. But from Lemma \ref{T906}(c) we have $A_0=A\cap S_0$ and
$A_1=A\cap S_1$, so $A_0,A_1\in V$. Thus either $S_0\cap\tau_G$ or $S_1\setminus\tau_G$ is a supported subset of
$\omega$, which is impossible. So $\sigma_G$ is unsupported and we are done. \kdok

In the same way one can prove the following statement concerning Black.

\begin{te}\bt\rm\label{T912}
Under the assumptions of Lemma \ref{T906}, Black wins the play $\langle p,p_0,i_0,$ $p_1,i_1,\dots\rangle$ in
the game

${\cal G}_1$ iff $\| \sigma\mbox{ is equal to }\check\omega \| >0$
iff $\exists q \leq p \;\; q \Vdash\tau = \check S_0$;

${\cal G}_2$ iff $\| \sigma$ is cofinite $\| >0$
iff $\exists q \leq p \;\; q\Vdash\tau = ^* \check S_0$;

${\cal G}_3$ iff $\| \sigma$ is supported $\| >0$
iff $\exists q \leq p \;\; q\Vdash ``\tau\cap \check S_0 \mbox{ or }
\check S_1\setminus\tau$ is supported";

${\cal G}_4$ iff $\| \sigma$ is infinite $\| >0$
iff $\exists q \leq p \;\; q \Vdash \tau \neq ^* \check S_1$.
\end{te}

\noindent Since for each sequence $\langle b_n \rangle$ in a c.B.a. ${\Bbb{B}}$
\begin{equation}\label{EQ900}
\textstyle
\bigwedge_{n\in\omega } b_n
\leq
\liminf b_n
\leq
\bigvee_{A\in [\omega ]^{\omega }}\bigwedge_{n\in A}b_n
\leq
\limsup b_n ,
\end{equation}
we have
%\begin{pro} \bt \rm \label{T900}
%Let the play  $\langle p,p_0,i_0,p_1,i_1,\dots\rangle$ satisfy the rules of Definition \ref{D900}. Then

%(a) White wins ${\cal G}_4$ $\Rightarrow$
%    White wins ${\cal G}_3$ $\Rightarrow$
%    White wins ${\cal G}_2$ $\Rightarrow$
%    White wins ${\cal G}_1$.

%(b) Black wins ${\cal G}_1$ $\Rightarrow$
%    Black wins ${\cal G}_2$ $\Rightarrow$
%    Black wins ${\cal G}_3$ $\Rightarrow$
%    Black wins ${\cal G}_4$.
%\end{pro}

\begin{pro} \bt \rm \label{T901}
Let ${\Bbb{B}}$ be a complete Boolean algebra. Then

(a) White has a w.s.\ in ${\cal G}_4$ $\Rightarrow$
    White has a w.s.\ in ${\cal G}_3$ $\Rightarrow$
    White has a w.s.\ in ${\cal G}_2$ $\Rightarrow$
    White has a w.s.\ in ${\cal G}_1$.

(b) Black has a w.s.\ in ${\cal G}_1$ $\Rightarrow$
    Black has a w.s.\ in ${\cal G}_2$ $\Rightarrow$
    Black has a w.s.\ in ${\cal G}_3$ $\Rightarrow$
    Black has a w.s.\ in ${\cal G}_4$.
\end{pro}

\section*{{\normalsize \bf 3. The game ${\cal G}_2$}}

\begin{te}\bt\rm\label{T909}
For each complete Boolean algebra ${\Bbb{B}}$ the following conditions are equivalent:

(a) ${\Bbb{B}}$ is not $(\omega,2)$-distributive;

(b) White has a winning strategy in the game ${\cal  G}_1 $;

(c) White has a winning strategy in the game ${\cal  G}_2  $.
\end{te}

\dok (a)$\Leftrightarrow$(b) is proved in \cite{Jech2} and (c)$\Rightarrow$(b) holds by Proposition \ref{T901}.
In order to prove (a)$\Rightarrow$(c) we suppose ${\Bbb{B}}$ is not $(\omega,2)$-distributive. Then
$p:=\|\exists x\subseteq\check\omega\;\;x\notin V\|>0$ and by The Maximum Principle there is a name $\pi\in
V^{{\Bbb{B}}}$ such that
\begin{equation}\label{EQ911}
p\Vdash\pi\subseteq\check\omega\;\land\;\pi\notin V.
\end{equation}
Clearly $\omega=A_0\cup A\cup A_p$, where $A_0=\{n\in\omega:\|\check n\in\pi\|\wedge p=0\}$,
$A=\{n\in\omega:\|\check n\in\pi\|\wedge p\in(0,p)_{{\Bbb{B}}}\}$ and $A_p=\{n\in\omega:\|\check n\in\pi\|\wedge
p=p\}$. We also have $A_0,A,A_p\in V$ and
\begin{equation}\label{EQ912}
p\Vdash\pi=(\pi\cap\check A)\cup\check A_p.
\end{equation}
Let $f:\omega\rightarrow A$ be a bijection belonging to $V$ and $\tau=\{\langle\check
n,\|f(n)\check{\enskip}\in\pi\|\land p\rangle:n\in\omega\}$. We prove
\begin{equation}\label{EQ913}
p\Vdash f[\tau]=\pi\cap\check A.
\end{equation}
Let $G$ be a ${\Bbb{B}}$-generic filter over $V$ containing $p$. If $n\in f[\tau_G]$ then $n=f(m)$ for some
$m\in\tau_G$, so $\|f(m)\check{\enskip}\in\pi\|\wedge p\in G$ which implies $\|f(m)\check{\enskip}\in\pi\|\in G$
and consequently $n\in\pi_G$. Clearly $n\in A$. Conversely, if $n\in\pi_G\cap A$, since $f$ is a surjection
there is $m\in\omega$ such that $n=f(m)$. Thus $f(m)\in\pi_G$ which implies $\|f(m)\check{\enskip}\in\pi\|\wedge
p\in G$ and hence $m\in\tau_G$ and $n\in f[\tau_G]$.

According to (\ref{EQ911}), (\ref{EQ912}) and (\ref{EQ913}) we have $p\Vdash\pi=f[\tau]\cup \check A_p\notin V$
so, since $A_p\in V$, we have $p\Vdash f[\tau]\notin V$ which implies $p\Vdash\tau\notin V$. Let
$p_n=\|f(n)\check{\enskip}\in\pi\|\wedge p$, $n\in\omega$. Then, by the construction, $p_n\in(0,p)_{{\Bbb{B}}}$
for all $n\in\omega$.

We define a strategy $\Sigma$ for White: at the beginning White plays $p$ and, in the $n$-th move, plays $p_n$.
Let us prove $\Sigma$ is a winning strategy for White in the game ${\cal  G}_2  $. Let $\langle
i_n:n\in\omega\rangle\in{}^{\omega }2$ be an arbitrary play of Black. According to Theorem \ref{T907} we prove
$p\Vdash\tau\neq^*\check S_0$. But this follows from $p\Vdash\tau\notin V$ and $S_0\in V$ and we are done.
\kdok

\section*{{\normalsize \bf 4. The game ${\cal G}_3$}}

Firstly we give some characterizations of complete Boolean algebras on which White has a winning strategy in
the game ${\cal G}_3$. To make the formulas more readable, we will write $w_\varphi$ for $w(\varphi)$. Also,
for $i:\omega \rightarrow 2$ we will denote $g^i=\{i\upharpoonright n:n\in\omega \}$, the corresponding
branch of the tree $^{<\omega }2$.

\begin{te} \bt \rm  \label{T902}
For a complete Boolean algebra ${\Bbb{B}}$ the following conditions are
equivalent:

(a) White has a winning strategy in the game ${\cal G}_3$  on ${\Bbb{B}}$;

(b) There are $p\in{\Bbb{B}}^+$ and $w:{}^{<\omega }2\rightarrow(0,p)_{{\Bbb{B}}}$
such that
\begin{equation}\label{EQ901}
\textstyle
\forall i:\omega \rightarrow 2\;
\bigvee_{A\in [\omega ]^{\omega }}\bigwedge_{n\in A}w_{i\upharpoonright n}^{i(n)}=0;
\end{equation}

(c) There are $p\in{\Bbb{B}}^+$ and $w:{}^{<\omega }2\rightarrow[0,p]_{{\Bbb{B}}}$
such that (\ref{EQ901}) holds.

(d) There are $p\in {\Bbb{B}}^+$ and $\rho \in V^{{\Bbb{B}}}$ such that

\parbox{11cm}
{\begin{eqnarray*}
& p\Vdash   \rho\subseteq(^{<\omega }2)\check{\enskip}\!\!\!\! &
\land\;
\forall\varphi\in(^{<\omega }2)\check{\enskip}\;
(\varphi ^{\smallfrown } \check{0} \in\rho\;\dot\lor\;\varphi ^{\smallfrown } \check{1}\in\rho) \\
&& \land\;
\forall i\in((^{\omega }2)^V)\check{\enskip}\; ( \rho\cap\check{g^i} \mbox{ is unsupported}).
\end{eqnarray*}}\hfill\parbox{1cm}{\begin{equation}\label{EQ905}\end{equation}}

(e) In some generic extension, $V_{{\Bbb{B}}}[G]$, there is a subset $R$ of the tree ${}^{< \omega }2$ containing
either $\varphi ^{\smallfrown } 0 $ or $\varphi ^{\smallfrown }1$,
for each $\varphi\in {}^{<\omega }2$, and having unsupported intersection with each branch of the
tree ${}^{< \omega }2$ belonging to $V$.
\end{te} \rm

\dok
(a)$\Rightarrow $(c). Let $\Sigma$ be a winning strategy for White.
$\Sigma$ is a function adjoining to each sequence of the form $\langle p,p_0,i_0,\dots,p_{n-1},i_{n-1}\rangle$,
where $p,p_0,\dots,p_{n-1}\in{\Bbb{B}}^+$ are obtained by $\Sigma$ and $i_0,i_1,\dots,i_{n-1}$
are arbitrary elements of $\{0,1\}$,
an element $p_n =\Sigma(\langle p,p_0,i_0,\dots,p_{n-1},i_{n-1}\rangle)$ of
$(0,p)_{{\Bbb{B}}}$ such that White playing in accordance with $\Sigma$ always wins. In general,
$\Sigma$ can be a multi-valued function, offering more ``good" moves for White, but according to The Axiom of
Choice, without loss of generality we suppose $\Sigma$ is a single-valued function, which is sufficient for the
following definition of $p$ and $w:{}^{<\omega }2\rightarrow[0,p]_{{\Bbb{B}}}$.

At the beginning $\Sigma$ gives $\Sigma(\emptyset)=p\in{\Bbb{B}}^+$ and, in the
first move, $\Sigma(\langle p\rangle)\in(0,p)_{{\Bbb{B}}}$. Let
$w_{\emptyset}=\Sigma(\langle p\rangle)$.

Let $\varphi\in {}^{n+1}2$ and let $w_{\varphi\upharpoonright k}$ be defined for $k\leq n$. Then we
define $w_{\varphi}=\Sigma(\langle p,w_{\varphi\upharpoonright 0},\varphi(0),\dots,w_{\varphi\upharpoonright
n},\varphi(n)\rangle)$.

In order to prove (\ref{EQ901}) we pick an $i:\omega \rightarrow 2$. Using induction it is easy
to show that in the match in
which Black plays $i(0),i(1),\dots,$ White,
following $\Sigma$ plays $p,w_{i\upharpoonright 0},w_{i\upharpoonright 1},\dots$ Thus,
since White wins ${\cal G}_3$, we have
$\bigvee_{A\in [\omega ]^{\omega }}\bigwedge_{n\in A}w_{i\upharpoonright n}^{i(n)}=0$
and (\ref{EQ901}) is proved.

(c)$\Rightarrow $(b). Let $p\in{\Bbb{B}}^+$ and $w:{}^{<\omega }2\rightarrow[0,p]_{{\Bbb{B}}}$ satisfy
(\ref{EQ901}). Suppose the set $S=\{\varphi\in {}^{<\omega }2:w_{\varphi}\in\{0,p\}\}$ is dense in the ordering
$\langle^{<\omega }2,\supseteq \rangle$. Using recursion we define $\varphi_k\in S$ for $k\in\omega $ as
follows. Firstly, we choose $\varphi_0\in S$ arbitrarily. Let $\varphi_k$ be defined and let $i_k\in 2$ satisfy
$i_k=0$ iff $w_{\varphi_k}=p$. Then  we choose $\varphi_{k+1}\in S$ such that $\varphi_k ^{\smallfrown }
i_k\subseteq\varphi_{k+1}$. Clearly the integers $n_k=\dom(\varphi_k)$, $k\in \omega$, form an increasing
sequence, so $i=\bigcup_{k\in\omega }\varphi_k:\omega \rightarrow 2$. Besides, $i\upharpoonright n_k=\varphi_k$
and $i(n_k)=i_k$. Consequently, for each $k\in\omega $ we have $w_{i\upharpoonright
n_k}^{i(n_k)}=w_{\varphi_k}^{i_k}=p$. Now $A_0 = \{ n_k : k\in \omega \} \in [\omega ]^{\omega }$ and
$\bigwedge_{n\in A_0}w_{i\upharpoonright n}^{i(n)}= p>0$. A contradiction to (\ref{EQ901}).

So there is $\psi\in {}^{<\omega }2$ such that $w_{\varphi}\in(0,p)_{{\Bbb{B}}}$, for all $\varphi\supseteq \psi$. Let
$m=\dom(\psi)$ and let $v_{\varphi}$ for $\varphi\in {}^{<\omega }2$ be defined by
$$v_{\varphi}=\left\{
\begin{array}{ll}
w_{\psi} & \mbox{if }|\varphi| < m,\\
w_{\psi {}^{\smallfrown }(\varphi\upharpoonright(\dom(\varphi)\setminus m))} & \mbox{if
}|\varphi| \geq m.
\end{array}\right.$$
Clearly $v: {}^{<\omega }2\rightarrow(0,p)_{{\Bbb{B}}}$ and we prove that $v$ satisfies (\ref{EQ901}). Let
$i:\omega \rightarrow 2$ and let $j=\psi ^{\smallfrown }(i\upharpoonright(\omega \setminus m))$. Then for $n\geq
m$ we have $v_{i\upharpoonright n}^{i(n)} =w_{\psi ^{\smallfrown } (i\upharpoonright(n\setminus m))}^{i(n)}
=w_{j\upharpoonright n}^{j(n)}$. Let $A\in [\omega ]^{\omega }$. Then $A\setminus m \in [\omega ]^{\omega }$
and, since $w$ satisfies (\ref{EQ901}), for the function $j$ defined above  we have $\bigwedge_{n\in A\setminus
m}w_{j\upharpoonright n}^{j(n)}=0$, that is $\bigwedge_{n\in A\setminus m}v_{i\upharpoonright n}^{i(n)}=0$,
which implies $\bigwedge_{n\in A}v_{i\upharpoonright n}^{i(n)}=0$ and (b) is proved.

(b)$\Rightarrow $(a). Assuming (b) we define a strategy $\Sigma$ for
White. Firstly White plays $p$ and $p_0=w_{\emptyset}$. In the $n$-th step, if $\varphi=\langle
i_0,\dots,i_{n-1}\rangle$ is the sequence of Black's previous moves, White plays $p_n=w_{\varphi}$. We prove
that $\Sigma$ is a winning strategy for White. Let $i:\omega \rightarrow 2$ code an arbitrary play of Black.
Since White follows $\Sigma$, in the $n$-th move White plays $p_n=w_{i\upharpoonright n}$, so according to
(\ref{EQ901}) we have
$\bigvee_{A\in [\omega ]^{\omega }}\bigwedge_{n\in A}p_n^{i_n}=0$
and White wins the game.

(b)$\Rightarrow $(d). Let $p\in{\Bbb{B}}^+$ and $w: {}^{<\omega }2\rightarrow(0,p)_{{\Bbb{B}}}$ be the
objects provided by (b). Let us define $v_{\emptyset}=p$ and, for $\varphi\in {}^{<\omega }2$ and $k\in 2$,
let $v_{\varphi ^{\smallfrown } k}=w_{\varphi}^k$. Then
$\rho=\{\langle\check\varphi,v_{\varphi}\rangle:\varphi\in {}^{<\omega }2\}$ is a name for a subset of
$^{<\omega }2$. If $i:\omega \rightarrow 2$, then $\sigma^i=\{\langle(i\upharpoonright
n)\check{\enskip},v_{i\upharpoonright n}\rangle:n\in\omega \}$ is a name for a subset of $g^i$ and, clearly,
\begin{equation}\label{EQ902}
1\Vdash\sigma^i=\rho\cap\check{g^i}.
\end{equation}
Let us prove
\begin{equation}\label{EQ903}
\forall i:\omega \rightarrow 2\;\;1\Vdash \rho\cap\check{g^i} \mbox{ is unsupported}.
\end{equation}
Let $i:\omega \rightarrow 2$. According to the definition of $v$, for $n\in\omega $ we have $w_{i\upharpoonright
n}^{i(n)}=v_{i\upharpoonright(n+1)}$ so, by (\ref{EQ901}), $\bigvee_{A\in[\omega ]^{\omega }}\bigwedge_{n\in
A}v_{i\upharpoonright(n+1)}=0$. By (\ref{EQ902}) we have
$v_{i\upharpoonright(n+1)}=\|(i\upharpoonright(n+1))\check{\enskip}\in\rho\cap\check{g^i}\|$ and we have $\|
\exists A \in (([\omega ]^{\omega })^V )\check{\enskip }\; \forall n\in A\;
(i\upharpoonright(n+1))\check{\enskip}\in\rho\cap\check{g^i}\|=0$ that is $\| \neg \exists B \in (([{}^{< \omega
}2 ]^{\omega })^V )\check{\enskip }\; B\subset \rho\cap\check{g^i}\|=1$ and (\ref{EQ903}) is proved. Now we
prove
\begin{equation}\label{EQ904}
\forall\varphi\in {}^{<\omega }2\;\; p\Vdash\check\varphi ^{\smallfrown }\check
0\in\rho \;\; \dot\lor \;\; \check\varphi ^{\smallfrown }\check 1\in\rho.
\end{equation}
If $p\in G$, where $G$ is a ${\Bbb{B}}$-generic filter over $V$,
then clearly $|G\cap\{w_{\varphi},p\setminus w_{\varphi}\}|=1$.
But $w_{\varphi}=w_{\varphi}^0=v_{\varphi ^{\smallfrown }0}=
\| \check\varphi ^{\smallfrown }\check 0\in\rho\|$ and
$p\setminus w_{\varphi}=w_{\varphi}^1=v_{\varphi{}^{\smallfrown }1}
=\|\check\varphi ^{\smallfrown }\check 1\in\rho\|$ and (\ref{EQ904}) is proved.

(d)$\Rightarrow $(c).
Let $p\in {\Bbb{B}}^+$ and $\rho \in V^{{\Bbb{B}}}$ satisfy (\ref{EQ905}).
In $V$ for each $\varphi\in {}^{<\omega }2$ we define
$w_{\varphi}=\|(\varphi ^{\smallfrown } 0)\check{\enskip}\in\rho\|\land p$ and check
condition (c). So for an arbitrary $i:\omega \rightarrow 2$ we prove
\begin{equation}\label{EQ906}
\textstyle
\bigvee_{A\in [\omega ]^{\omega }}\bigwedge _{n\in A}w_{i\upharpoonright n}^{i(n)}=0.
\end{equation}
According to (\ref{EQ905}) for each $n\in\omega $ we have
$p\Vdash((i\upharpoonright n)^{\smallfrown } 0)\check{\enskip}\in\rho \; \dot\lor \;
((i\upharpoonright n)^{\smallfrown }1)\check{\enskip}\in\rho$, that is $p\leq a_0\lor a_1$ and $p\land a_0\land
a_1=0$, where $a_k=\|((i\upharpoonright n)^{\smallfrown } k)\check{\enskip}\in\rho\|$, $k\in\{0,1\}$,
which clearly implies $p\land a_0'=p\land a_1$, i.e.
\begin{equation}\label{EQ907}
p\land\|((i\upharpoonright n)^{\smallfrown } 0)\check{\enskip}\in\rho\|'=p
\land\|((i\upharpoonright n)^{\smallfrown } 1)\check{\enskip}\in\rho\|.
\end{equation}
Let us prove
\begin{equation}\label{EQ908}
w_{i\upharpoonright n}^{i(n)}=\|(i\upharpoonright(n+1))\check{\enskip}\in\rho\|\land p.
\end{equation}
If $i(n)=0$, then $w_{i\upharpoonright n}^{i(n)}=
\|((i\upharpoonright n) ^{\smallfrown } 0)\check{\enskip}\in\rho\|\land p=\|((i\upharpoonright n)^{\smallfrown }
i(n))\check{\enskip}\in\rho\|\land p$ and (\ref{EQ908}) holds.
If $i(n)=1$, then according to (\ref{EQ907}) $w_{i\upharpoonright
n}^{i(n)}=p\setminus w_{i\upharpoonright n}=
p\land\|((i\upharpoonright n)^{\smallfrown } 0)\check{\enskip}\in\rho\|'=
p\land\|((i\upharpoonright n)^{\smallfrown }1)\check{\enskip}\in\rho\|=
p\land\|((i\upharpoonright n)^{\smallfrown } i(n))\check{\enskip}\in\rho\|$ and (\ref{EQ908}) holds again.

Now $ \bigvee_{A\in [\omega ]^{\omega }}\bigwedge _{n\in A} w_{i\upharpoonright n}^{i(n)} =p\land \| \exists A
\in (([\omega ]^{\omega })^V )\check{\enskip } \; \forall  n\in A\; \check{i}\upharpoonright(n+1)\in\rho\|
=p\land\| \rho\cap\check{g^i} \mbox{ is supported} \|=0$, since by (\ref{EQ905}) $p \leq \| \rho\cap\check{g^i}
\mbox{ is unsupported}\|$. Thus (\ref{EQ906}) is
proved.

(d)$\Rightarrow$(e) is obvious and (e)$\Rightarrow$(d) follows from The Maximum Principle.
\kdok

Concerning condition (e) of the previous theorem we note that in \cite{KuSo} the following characterization
is obtained.

\begin{te} \bt \rm  \label{T913}
White has a winning strategy in the game ${\cal G}_4$  on a c.B.a.\ ${\Bbb{B}}$
if and only if in some generic extension, $V_{{\Bbb{B}}}[G]$, there is a subset
$R$ of the tree ${}^{< \omega }2$ containing
either $\varphi ^{\smallfrown }0$ or $\varphi ^{\smallfrown }1$,
for each $\varphi\in {}^{<\omega }2$,
and having finite intersection with each branch of the
tree ${}^{< \omega }2$ belonging to $V$.
\end{te} \rm

\begin{te}\bt\rm\label{T904}
Let $\Bbb B$ be a complete Boolean algebra. If forcing by ${\Bbb{B}}$ produces an independent real in some
generic extension, then White has a winning strategy in the game ${\cal G}_3$ played on ${\Bbb{B}}$.
\end{te}

\dok
Let $p=\|\exists x\subseteq\check\omega\;\; x\mbox{ is independent}\|>0$. Then, by The Maximum Principle
there is a name $\tau\in V^{\Bbb{B}}$ such that
\begin{equation}\label{EQ909}
p\Vdash\tau\subseteq\check\omega\;\land\;\forall
A\in(([\omega]^{\omega})^V)\check{\enskip}\;\;(|A\cap\tau|=\check\omega\;\land\;|A\setminus\tau|=\check\omega).
\end{equation}
Let us prove that $K=\{n\in\omega:\|\check n\in\tau\|\land p\in\{0,p\}\}$ is a finite set. Clearly $K=K_0\cup
K_p$, where $K_0=\{n\in\omega:p\Vdash\check n\notin\tau\}$ and $K_p=\{n\in\omega:p\Vdash\check n\in\tau\}$.
Since $p\Vdash \check K_0\subseteq\check\omega\setminus\tau\;\land\;\check K_p\subseteq\tau$, according to
(\ref{EQ909}) the sets $K_0$ and $K_p$ are finite, thus $|K|<\omega$.

Let $q\in(0,p)_{{\Bbb{B}}}$ and let $p_n$, $n\in\omega$, be defined by
$$p_n=\left\{\begin{array}{ll}
q & \mbox{if }n\in K, \\
\|\check n\in\tau\|\land p & \mbox{if }n\in\omega\setminus K.
\end{array}\right.$$
Then for $\tau_1=\{\langle\check n,p_n\rangle:n\in\omega\}$ we have $p\Vdash\tau_1=^*\tau$ so according to
(\ref{EQ909})
\begin{equation}\label{EQ910}
p\Vdash\tau_1\subseteq\check\omega\;\land\;\forall
A\in(([\omega]^{\omega})^V)\check{\enskip}\;\;(|A\cap\tau_1|=\check\omega\;\land\;|A\setminus\tau_1|=\check\omega).
\end{equation}
Then $p_n=\|\check n\in\tau_1\|\in(0,p)_{{\Bbb{B}}}$ and we define a strategy $\Sigma$ for White: at the
beginning White plays $p$ and, in the $n$-th move, White plays $p_n$.

We prove $\Sigma$ is a winning strategy for White. Let $\langle p,p_0,i_0,p_1,i_1,\dots\rangle$ be an arbitrary
play in which White follows $\Sigma$ and let $S_k=\{n\in\omega:i_n=k\}$, for $k\in\{0,1\}$. Suppose
$q=\bigvee_{A\in[\omega]^{\omega}}\bigwedge_{n\in A}p_n^{i_n}>0$. Now $q\leq p$ and $q=
\bigvee_{A\in[\omega]^{\omega}}(\bigwedge_{n\in A\cap S_0}\|\check n\in\tau_1\|\; \land\;\bigwedge_{n\in A\cap
S_1}(p\land\|\check n\notin\tau_1\|) =p\land\bigvee_{A\in[\omega]^{\omega}}\|\check A\cap\check
S_0\subseteq\tau_1\; \land\;\check A\cap\check S_1\subseteq \check \omega \setminus \tau_1\| \leq\|\exists
A\in(([\omega]^{\omega})^V)\check{\enskip}\;(\check A\cap\check S_0\subseteq\tau_1\; \land\;\check A\cap\check
S_1\subseteq \check\omega \setminus\tau_1)\|$.

Let $G$ be a ${\Bbb{B}}$-generic filter over $V$ containing $q$. Then there is  $A\in[\omega]^{\omega}\cap V$
such that $A\cap S_0\subseteq(\tau_1)_G$ and $A\cap S_1\subseteq\omega\setminus(\tau_1)_G$. But one of the sets
$A\cap S_0$ and $A\cap S_1$ must be infinite and, since $p\in G$, according to (\ref{EQ910}), it must be
split by $(\tau_1)_G$. A contradiction. Thus $q=0$ and White wins the game.
\kdok

\begin{te}\bt\rm\label{T908}
Let ${\Bbb{B}}$ be an $(\omega,2)$-distributive complete Boolean algebra. Then

(a) If $\langle p,p_0,i_0,p_1,i_1,\dots\rangle$ is a play satisfying the rules given in Definition \ref{D900},
then Black wins the game ${\cal G}_3$ iff Black wins the game ${\cal G}_4$.

(b) Black has a winning strategy in the game ${\cal G}_3$ iff Black has a winning strategy in the game ${\cal G}_4$.
\end{te}

\dok (a) The implication ``$\Rightarrow$" follows from the proof of Proposition \ref{T901}(b). For the proof
of ``$\Leftarrow$" suppose Black wins the play $\langle p,p_0,i_0,p_1,i_1,\dots\rangle$ in the game ${\cal
G}_4$. Then, by Theorem \ref{T912} there exists $q\in{\Bbb{B}}^+$ such that $q\Vdash$``$\sigma$ is
infinite". Since the algebra ${\Bbb{B}}$ is $(\omega,2)$-distributive we have $1\Vdash\sigma\in V$, thus
$q\Vdash\sigma\in(([\omega]^{\omega})^V)\check{\enskip}$ and hence $\neg 1\Vdash``\sigma$ is not supported"
so, by Theorem \ref{T912}, Black wins ${\cal G}_3$.

(b) follows from (a).
\kdok

\section*{{\normalsize \bf 5. Indeterminacy, problems}}

\begin{te}\bt\rm\label{T911}
$\diamondsuit$ implies the existence of a Suslin algebra on which the games
${\cal  G}_1$,${\cal  G}_2$, ${\cal  G}_3$ and ${\cal  G}_4 $ are undetermined.
\end{te}

\dok
Let ${\Bbb{B}}$ be the Suslin algebra mentioned in (c) of Theorem \ref{T910}. According to Proposition
\ref{T901}(b) and since Black does not have a winning strategy in the game ${\cal  G}_4    $,
Black does not have a winning
strategy in the games ${\cal  G}_1 ,{\cal  G}_2  ,{\cal  G}_3   $ as well. On the other hand,
since the algebra ${\Bbb{B}}$ is
$(\omega,2)$-distributive, White does not have a winning strategy in the game ${\cal  G}_1 $ and, by Proposition
\ref{T901}(a), White does not have a winning strategy in the games ${\cal  G}_2  ,{\cal  G}_3   ,{\cal  G}_4    $
played on ${\Bbb{B}}$.
\kdok

\begin{pb}\bt\rm\label{P901}
According to Theorem \ref{T904}, Proposition \ref{T901} and Theorem \ref{T909} for each complete Boolean
algebra ${\Bbb{B}}$ we have:
\begin{center}
${\Bbb{B}}$ is $\omega$-independent $\Rightarrow$ White has a winning strategy in ${\cal  G}_3   $ $\Rightarrow$
${\Bbb{B}}$ is not $(\omega,2)$-distributive.
\end{center}
Can one of the implications be reversed?
\end{pb}

\begin{pb}\bt\rm\label{P902}
According to Proposition \ref{T901}(b), for each complete Boolean algebra ${\Bbb{B}}$ we have:
\begin{center}
Black has a winning strategy in ${\cal  G}_1 $ $\Rightarrow$ Black has a winning strategy
in ${\cal  G}_2  $ $\Rightarrow$ Black
has a winning strategy in ${\cal  G}_3   $.
\end{center}
Can some of the implications be reversed?
\end{pb}

We note that the third implication from Proposition \ref{T901}(b) can not be replaced by the
equivalence, since if ${\Bbb{B}}$ is the Cohen or the random algebra, then Black has a winning strategy in the
game ${\cal  G}_4$ (Theorem \ref{T910}(b)) while Black does not have a winning strategy in the game ${\cal  G}_3$,
because White has one (the Cohen and the random forcing produce independent reals and Theorem \ref{T904} holds).

\end{document}